\documentclass[11pt]{amsart}
 
\usepackage{amsfonts}
\usepackage{amssymb}
\usepackage{amscd}
\input xy
\xyoption{all}
 
\newcommand{\marginlabel}[1]%
  {\mbox{}\marginpar{\raggedleft\hspace{0pt}\bfseries\sf#1}}

\theoremstyle{plain}
\newtheorem{thm}{Theorem}[section]
\newtheorem{lem}[thm]{Lemma}
\newtheorem{prop}[thm]{Proposition}
\newtheorem{cor}[thm]{Corollary}

\theoremstyle{definition}

\theoremstyle{remark}
\newtheorem{eg}[thm]{Example}

\newtheorem{rmk}[thm]{Remark}

\numberwithin{equation}{section}

\def\N{{\mathbb N}}

\def\A{{\mathbb A}}
\def\R{{\mathbb R}}

\def\P{{\mathbb P}}

\def\H{\mathcal{H}}

\def\O{\mathcal{O}}

\def\J{\mathcal{J}}

\def\PP{\mathcal{P}}
\def\QQ{\mathcal{Q}}

\def\a{\mathfrak{a}}
\def\bb{\mathfrak{b}}
\def\cc{\mathfrak{c}}                           
\def\dd{\mathfrak{d}}                           
\def\mm{\mathfrak{m}}
\def\DDD{\mathfrak{D}}

\def\MM{\mathfrak{M}}                          

\def\cc{{\bf c}}

\def\g{\gamma}
\def\d{\delta}
\def\f{\phi}

\def\e{\eta}

\def\n{\nu}
\def\m{\mu}

\def\p{\pi}

\def\D{\Delta}

\def\LL{\Lambda}
\def\S{\Sigma}

\def\o{\circ}
\def\ov{\overline}
\def\rat{\dashrightarrow}                      

\def\inj{\hookrightarrow}

\def\.{\cdot}
\def\({\Big{(}}
\def\){\Big{)}}
\def\^{\widehat}
\def\-{}

\DeclareMathOperator{\codim} {codim}

\DeclareMathOperator{\Proj} {Proj}
\DeclareMathOperator{\Spec} {Spec}

\DeclareMathOperator{\Aut} {Aut}
\DeclareMathOperator{\Bir} {Bir}

\DeclareMathOperator{\Bl} {Bl}

\DeclareMathOperator{\ord} {ord}

\DeclareMathOperator{\Supp} {Supp}

\DeclareMathOperator{\ini} {in}

\DeclareMathOperator{\lc} {lc}

\DeclareMathOperator{\pr} {pr}

\begin{document}

\title{Bounds for log canonical thresholds with applications to
  birational rigidity}

\author[T. de Fernex]{Tommaso de Fernex}
\address{Department of Mathematics, University of Michigan,
East Hall, 525 East University Avenue, Ann Arbor, MI 48109-1109, USA}
\email{{\tt defernex@umich.edu}}

\author[L. Ein]{Lawrence Ein}
\address{Department of Mathematics, University of Illinois at Chicago,
851 S. Morgan St., M/C. 249, Chicago, IL 60607-7045, USA}
\email{{\tt ein@math.uic.edu}}

\author[M. Musta\c{t}\v{a}]{Mircea~Musta\c{t}\v{a}}
\address{Department of Mathematics, Harvard University,
1 Oxford Street, Cambridge, MA 02138, USA}
\email{{\tt mirceamustata@yahoo.com}}

\maketitle
\markboth{T. de Fernex, L.~Ein and M. Musta\c{t}\v{a}}
{\bf BOUNDS FOR LOG CANONICAL THRESHOLDS}

\section*{Introduction}

Let $X$ be a smooth algebraic variety, defined over an
algebraically closed field of characteristic zero, 
and let $V \subset X$ be a proper closed subscheme.
Our main goal in this paper is to study an invariant
of the pair $(X,V)$, called
the log canonical threshold of $X$ along $V$, and
denoted by $\lc(X,V)$.
Interest in bounds for log canonical thresholds is motivated by 
techniques that have recently been developed in 
higher dimensional birational geometry. 
In this paper, we study this invariant using intersection
theory, degeneration techniques and jet schemes.

A natural question is how does this invariant behave
under basic operations such as restrictions
and projections. Restriction properties have been 
extensively studied in recent years, leading to important
results and conjectures. In the first section of this paper,
we investigate the behavior
under projections, and we prove the
following result (see Theorem~\ref{thm1} for a more precise statement):

\begin{thm}\label{thm1-intro}
With the above notation, suppose that $V$ is Cohen-Macaulay,
of pure codimension $k$,
and let $f : X \to Y$ be a proper, dominant, smooth morphism of 
 relative dimension $k-1$, with $Y$ smooth.
If $f|_V$ is finite, then 
$$
\lc(Y,f_*[V]) \leq \frac{k! \. \lc(X,V)^k}{k^k},
$$
and the inequality is strict if $k\geq 2$.
Moreover, if $V$ is locally complete intersection,
then 
$$
\lc(Y,f_*[V]) \leq \frac{\lc(X,V)^k}{k^k}.
$$
\end{thm}

Examples show that these bounds are sharp.
The proof of the above theorem is based on a general
inequality relating the log canonical threshold of
a fractional ideal of the form $h^{-b}\. \a$,
and the colength of $\a$. Here $\a$ is a zero dimensional
ideal in the local ring of $X$ at some (not necessarily closed) point,
$b\in{\mathbb Q}_+$, and $h$ is the equation of a smooth divisor.
We prove this inequality in the second section (see Theorem~\ref{l(a)-e(a)}),
using a degeneration to monomial ideals. It generalizes 
a result from \cite{DEM}, which was the case $b=0$.

In the third section, we give lower bounds for the log canonical threshold
of affine subschemes defined by homogeneous equations of the same degree.
We prove the following

\begin{thm}
Let $V\subset X=\A^n$ be a subscheme defined by homogeneous equations of 
degree $d$. Let $c=\lc(\A^n, V)$, and let $Z$ be the non log terminal locus of
$(\A^n, c\. V)$. If $e=\codim(Z,\A^n)$, then
$$
\lc(\A^n,V) \ge \frac{e}d.
$$
Moreover, we have equality if and only if the following holds:
$Z$ is a linear subspace, and if $\pi : \A^n\longrightarrow\A^n/Z$
is the projection, then there is a subscheme $V'\subset\A^n/Z$
such that $V=\pi^{-1}(V')$, $\lc(\A^n/Z,V')=e/d$, and the
non log terminal locus of $(\A^n/Z,(e/d)\. V')$ is the origin.
\end{thm}

The proof of this result is based on the characterization of the log canonical
threshold via jet schemes from~\cite{Mu2}.
In the particular case when $V$ is the affine cone over a 
projective hypersurface with isolated singularities,
the second assertion in the above result proves a conjecture of Cheltsov and Park
from~\cite{CP}.

In the last section we apply the above bounds in the context
of birational geometry. In their influential paper \cite{IM},
Iskovskikh and Manin proved that a smooth quartic threefold is
what is called nowadays 
birationally superrigid; in particular, every birational 
automorphism is regular, and the variety is not rational. There has
been a lot of work to extend this result to other Fano varieties of index one,
in particular to smooth hypersurfaces of degree $N$ in $\P^N$, for $N>4$.
The case $N=5$ was done by Pukhlikov in \cite{Pu2}, and the cases
$N=6,7,8$ were proven by Cheltsov in \cite{Ch2}.
Moreover,
Pukhlikov showed in \cite{Pu5} that a 
general hypersurface as above is birationally superrigid, for every $N>4$. 
We use our results to give an easy and uniform proof of birational
superrigidity for arbitrary smooth hypersurfaces of degree $N$ in $\P^N$
when $N$ is small.

\begin{thm}\label{thm3_introd}
If $X\subset{\mathbb P}^N$ is a smooth hypersurface of degree $N$, and if
$4\leq N\leq 12$, then $X$ is birationally superrigid. 
\end{thm} 

Based on previous ideas of Corti, Pukhlikov proposed in \cite{Pu1}
a proof of the birational rigidity of every smooth hypersurface 
of degree $N$ in $\P^N$, for $N\geq 6$.
Unfortunately, at the moment there is a gap in his arguments (see
Remark~\ref{gap} below).
Despite this gap, the proof proposed in \cite{Pu1} contains many
remarkable ideas, and
it seems likely that a complete proof could be obtained in the
future along those lines. In fact, the outline of the proof of
Theorem~\ref{thm3_introd} follows his method,
and our contribution is mainly to simplifying and solidifying his argument.

\subsection*{Acknowledgements}
We are grateful to Steve Kleiman and Rob Lazarsfeld for
useful discussions.
Research of the first author was partially supported
by MURST of Italian Government, National Research Project 
(Cofin 2000) ``Geometry of Algebraic Varieties''.
Research of the second author was partially supported by 
NSF~Grant DMS~02-00278. The third author
served as a Clay Mathematics Institute Long-Term Prize Fellow
while this research has been done.

\section{Singularities of log pairs under projections}

Let $X$ be a smooth algebraic variety, defined over an algebraically 
closed field of characteristic zero, and let $V \subset X$ 
be a proper subscheme.
For any rational number $c > 0$, we can consider the pair $(X,c\. V)$.
The usual definitions in the theory of singularities of pairs,
for which we refer to \cite{Ko}, extend to this context.
In particular, we say that 
an irreducible subvariety $C \subset X$ is 
a center of non log canonicity
(resp. non log terminality, non canonicity, non terminality)
for $(X,c\. V)$ if there is at least one divisorial valuation of $K(X)$, 
with center $C$ on $X$, whose discrepancy along $(X,c\.V)$ is 
$<-1$ (resp. $\le -1$, $<0$, $\le 0$). We will denote by 
$\lc(X,V)$ the log canonical threshold of the pair $(X,V)$, i.e.,
the largest $c$ such that $(X,c\. V)$ is log canonical.
We will occasionally consider also pairs of the form $(X,c_1\. V_1-c_2\.V_2)$,
where $V_1$, $V_2\subset X$ are proper subschemes of $X$.
The definition of (log) terminal and canonical pairs extends in an obvious
way to this setting.

We fix now the set-up for this section.
Let $f : X \to Y$ be a smooth and proper morphism onto 
a smooth algebraic variety $Y$. We assume that $V\subset X$ is 
a pure dimensional, Cohen-Macaulay closed subscheme, such that
$\dim V = \dim Y - 1$, and such that
the restriction of $f$ to $V$ is finite.
If $[V]$ denotes the cycle associated to $V$, then its
push-forward 
$f_*[V]$ determines an effective Cartier divisor on $Y$. We set 
$\codim(V,X)=k$.

\begin{thm}\label{thm1}
With the above notation, let  
$C \subset X$ be an irreducible center of
non log terminality for $(X,c\. V)$, for some $c>0$. Then
$f(C)$ is a center of non log terminality
(even non log canonicity, if $k\geq 2$)
for the pair
\begin{equation}\label{gen_formula}
\( Y, \frac{k! \. c^k}{k^k} \. f_*[V] \).
\end{equation}
Moreover, if $V$ is locally complete intersection
(l.c.i. for short) then 
$f(C)$ is a center of non log terminality for the pair
\begin{equation}\label{lci_formula}
\( Y, \frac{c^k}{k^k} \. f_*[V] \).
\end{equation}
\end{thm}

\begin{eg}
Let $k$ and $n$ be two positive integers with $n > k$, and let 
$R = K[x_k,\dots,x_n]$. We take
$X = \P^{k-1}_R = \Proj R[x_0,\dots,x_{k-1}]$, 
$Y = \Spec R$, and let $f$ be the natural
projection from $X$ to $Y$. 
For any $t>0$, let $V_t$ be the subscheme of $X$ defined by the
homogeneous ideal 
$(x_1,\dots,x_k)^t$. Note that 
$\lc(X,V_t) = k/t$, and that if $c=k/t$, then $V_1$ is a center of
non log terminality for $(X,c\. V_t)$.
Since $l(\O_{V_t,V_1})=\binom{k+t-1}{k}$, we see that
$$\lim_{t\to\infty}\frac{k! \. c^k/k^k}{\lc(Y,f_*[V_t])}
=\lim_{t\to\infty}\frac{t(t+1)\ldots(t+k-1)}{t^k}=1,$$ 
so the bound in (\ref{gen_formula}) is sharp
(at least asymptotically).

To prove sharpness in the l.c.i. case, let
$W_t\subset X$ be the complete intersection subscheme defined
by $(x_1^t,\dots,x_k^t)$. This time $l(\O_{W_t,W_1}) = t^k$, and
$\lc(Y,f_*[W_t]) = 1/t^k = \lc(X,W_t)^k/k^k$.
\end{eg}

\begin{proof}[Proof of Theorem~\ref{thm1}]
By hypothesis, there is a proper
birational morphism $\n : W \to X$, where $W$ can be chosen
to be smooth, and a smooth irreducible divisor $E$ on $W$, 
such that $\n(E) = C$, and
such that 
the discrepancy of $(X,c\. V)$ at $E$ is 
\begin{equation}\label{eq1}
a_E(X,c\. V) \le -1.
\end{equation}
The surjection $f$ induces an inclusion of function fields
 $f^* : K(Y) \inj K(X)$.
Let $R_E:=\O_{W,E}\subset K(X)$ be the discrete valuation ring associated 
to the valuation along $E$, and let $R = (f^*)^{-1}R_E$. 
Note that $R$ is a non-trivial discrete valuation ring.

\begin{lem}
$R$ corresponds to a divisorial valuation.
\end{lem}

\begin{proof}
It is enough to show that the transcendence degree of the
residue field of $R$ over the ground field is
$\dim Y-1$ (see~\cite{KM}, Lemma~2.45).
This follows from~\cite{ZS}, VI.6, Corollary~1.
\end{proof}

The lemma implies that there is a proper birational
morphism $\g : Y' \to Y$ and an irreducible divisor $G$ on $Y'$ such that 
$R=\O_{Y',G}$.
By Hironaka's theorem, we may assume that both $Y'$ and $G$ are smooth,
and moreover, that the union between $G$ and the exceptional locus of $\g$
has simple normal crossings.
Since the center of $R_E$ on $X$ is $C$, we deduce that $R$ has center
$f(C)$ on $Y$, so $\g(G) = f(C)$. 

Consider the fibered product
$X' = Y' \times_Y X$.
We may clearly 
assume that $\n$ factors through the natural map $\f : X' \to X$.
Therefore we have the following commutative diagram:
$$
\xymatrix{
W \ar[r]^{\e} & X' \ar[d]_g \ar[r]^{\f} & X \ar[d]^f \\
&Y' \ar[r]^{\g} & Y,
}
$$
where $\phi\circ\eta=\nu$.
Note that $X'$ is a smooth variety,
 $g$ is a smooth, proper morphism, and $\e$ and $\f$
are proper, birational morphisms. Let $V' = \f^{-1}(V)$
be the scheme theoretic inverse image of $V$ in $X'$, 
i.e.,  the subscheme of $X'$ defined by the ideal sheaf $I_V \. \O_{X'}$.
 
\begin{lem}\label{lem1}
$V'$ is pure dimensional, $\codim(V',X')=k$, and 
$\f^*[V]$ is the class of $[V']$.
Moreover, if $V$ is l.c.i., then so is $V'$.
\end{lem}

\begin{proof}
Note that both 
$\gamma$ and $\phi$ are l.c.i. morphisms, because they are morphisms between
smooth varieties.  
The pull-back in the statement is the pull-back by such a 
 morphism (see \cite{Fulton}, Section 6.6). Recall how this 
is defined. We factor $\gamma$ as $\gamma_1\circ\gamma_2$, where
$\gamma_1 : Y'\times Y\longrightarrow Y$ is the projection, and 
$\gamma_2 : Y'\hookrightarrow Y'\times Y$ is the graph of $\gamma$. 
By pulling-back, we get a corresponding decomposition $\phi=\phi_1\circ\phi_2$,
with $\phi_1$ smooth, and $\phi_2 : X'\hookrightarrow Y'\times X$ 
a regular embedding of 
codimension $\dim Y'$.
Then $\phi^*[V]=\phi_2^!([Y'\times V])$.

Since $f|_V$ is finite and $V' = Y'\times_Y V$, $g|_{V'}$ is also 
finite. Moreover, since $g(V')$ is a proper subset of $Y'$,
we see that $\dim V' \le \dim Y'-1$. 
On the other hand, $V'$ is locally cut in $Y'\times V$
by $\dim\,Y'$ equations, so that every irreducible component of $V'$
has dimension at least $\dim\,V$. Therefore $V'$ is pure dimensional, and 
$\dim\,V'=\dim\,V$.

Since $Y'\times V$ is Cohen-Macaulay, this also
implies that $\phi_2^!([Y'\times V])$ is equal to the class of
$[V']$, by Proposition~7.1 in
\cite{Fulton}. This proves the first assertion.
Moreover, if $V$ is l.c.i.,
then it is locally defined in $X$ by $k$ equations. The same 
is true for $V'$, hence
$V'$ is l.c.i., too.
\end{proof}

We will use the following notation for multiplicities.
Suppose that $W$ is an irreducible subvariety of a variety $Z$.
Then the multiplicity of $Z$ along $W$ is denoted by $e_WZ$
(we refer to \cite{Fulton}, Section 4.3, for definition
and basic properties). If $\alpha=\sum_in_i[T_i]$ is a pure
dimensional cycle on $Z$, then $e_W\alpha:=\sum_in_ie_WT_i$
(if $W\not\subseteq T_i$, then we put $e_WT_i=0$).
Note that if $W$ is a prime divisor, and if $D$
is an effective Cartier divisor on $Z$, then we have
$e_W[D]={\rm ord}_W(D)$, where $[D]$ is the cycle associated to $D$,
and ${\rm ord}_W(D)$ is the coefficient of $W$ in $[D]$. As we work
on smooth varieties, from now on we will identify $D$ with $[D]$.

Let $F = \e(E)$. Note that by construction, we have $g(F)=G$.
Since $F\subseteq V'$, and $g|_{V'}$ is finite, and $\dim\,G=\dim\,V'$,
it follows that $F$ is an irreducible component of $V'$,
hence $\codim(F, X')=k$.
We set $a = e_F(K_{X'/X})$.

To simplify the statements, we put
$$
\d = 
\begin{cases}
1 &\text{if $V$ is a l.c.i.,} \\
k! &\text{otherwise.}
\end{cases}
$$

\begin{lem}\label{lem2}
We have
$$
\ord_G(\g^* f_*[V]) \geq \frac{(a + 1)k^k}{\d c^k},
$$
and the inequality is strict in the case $\delta=k!$, if $k\geq 2$.
\end{lem}

\begin{proof}
Since $\f$ and $\g$ are l.c.i. morphisms
of the same relative dimension, it follows from
\cite{Fulton}, Example~17.4.1, and Lemma~\ref{lem1} that
$g_*[V']$ and $\g^*f_*[V]$ are linearly equivalent, as divisors on $Y'$.
As the two divisors are equal outside the exceptional locus of $\g$,
we deduce from the Negativity Lemma 
(see \cite{KM}, Lemma 3.39) that also 
their $\g$-exceptional components must coincide. 
This gives  $g_*[V'] = \g^*f_* [V]$.

In particular,
${\rm ord}_G(\g^*f_*[V])$
is greater or equal to the coefficient of $F$ in $[V']$.
Lemma~\ref{lem1} implies
$$
\ord_G(\g^* f_*[V]) \geq l(\O_{V',F}),
$$
so that it is enough to show that
\begin{equation}\label{lem2-eq}
l(\O_{V',F}) \geq \frac{(a + 1)k^k}{\d c^k},
\end{equation}
and that the inequality is strict in the case $\delta=k!$, if $k\geq 2$.

By replacing $W$ with a higher model, we may clearly
assume that $\n^{-1}(V)$ is an effective divisor on $W$.
If $I_V\subseteq\O_X$ is the ideal defining $V$, then
we put $\ord_E(I_V):=\ord_E\n^{-1}(V)$. It follows from~(\ref{eq1}) that we have 
$$
-1\geq \ord_E(K_{W/X}) - c\.\ord_E(I_V) = \ord_E(K_{W/X'}) - 
(c\.\ord_E(I_{V'})-\ord_E(K_{X'/X})).
$$
Therefore $F$ is a center of non log terminality for the 
pair $(X',c\. V' - K_{X'/X})$. Since $g(F)=G$ is a divisor on $Y'$,
it follows that $F$ can not be contained in the intersection of two distinct
$\phi$-exceptional divisors. Hence
the support
of $K_{X'/X}$  is smooth at the generic point of $F$.
Then~(\ref{lem2-eq}) follows from
Theorem~\ref{l(a)-e(a)} below (note that the 
length of a complete intersection 
ideal coincides with its Samuel multiplicity).
\end{proof}

We continue the proof of Theorem~\ref{thm1}.
Note that $\ord_G K_{Y'/Y} \leq e_F(g^* K_{Y'/Y})$. 
Since $K_{X'/X} = g^* K_{Y'/Y}$ (see~\cite{Hartshorne},
Proposition~II~8.10), we deduce 
$$
\ord_G K_{Y'/Y} \leq a.
$$
In conjunction with Lemma~\ref{lem2}, this gives 
$$
a_G\left(Y,\frac{\delta c^k}{k^k}f_*[V]\right)=
\ord_G \( K_{Y'/Y} -  \frac{\d c^k}{k^k}\. \g^* f_*[V] \) \leq -1. 
$$
Moreover, this inequality is strict in the case when $\d=k!$, if $k\geq 2$.
This completes the proof of Theorem~\ref{thm1}.
\end{proof}

\begin{rmk}
We refer to~\cite{Pu1} for a result on the canonical threshold of complete
intersection subschemes of codimension 2, via generic projection.
\end{rmk}

\bigskip

\section{Multiplicities of fractional ideals}

In this section we extend some of the results
of~\cite{DEM}, as needed in the proof of Theorem~\ref{thm1}. 
More precisely, we consider the following set-up.
Let $X$ be a smooth variety, $V\subset X$ a closed subscheme, and
let $Z$ be an irreducible component of $V$. We denote by
$n$ the codimension of $Z$ in $X$, and by
$\a \subset\O_{X,Z}$ the image of the ideal defining $V$. 
Let $H \subset X$ be a prime 
divisor containing $Z$, such that $H$
is smooth at the generic point of $Z$.
We consider the pair
$$
(X, V-b\cdot H),
$$
for a given $b\in{\mathbb Q}_+$. 

\begin{thm}\label{l(a)-e(a)}
With the above notation, suppose that for some $\mu\in{\mathbb Q}_+^*$,
$(X,\frac{1}{\mu}(V-b\cdot H))$ is not log terminal
at the generic point of $Z$. Then
\begin{equation}\label{l(a)}
l(\O_{X,Z}/\a)\geq\frac{n^n \m^{n-1}(\m + b)}{n!},
\end{equation}
and the inequality is strict if $n\geq 2$.
Moreover, if $e(\a)$ denotes the Samuel multiplicity of $\O_{X,Z}$
along $\a$, then
\begin{equation}\label{e(a)}
e(\a) \ge n^n \m^{n-1}(\m+b).
\end{equation}
\end{thm}

\begin{rmk}
For $n=2$, inequality~(\ref{e(a)}) gives a result of 
Corti from~\cite{Co2}. On the other hand, if $b=0$, then the
statement reduces to Theorems~1.1 and~1.2 in \cite{DEM}.
\end{rmk}

\begin{proof}[Proof of Theorem~\ref{l(a)-e(a)}]
We see that~(\ref{l(a)}) implies~(\ref{e(a)}) as follows. If we apply the
first formula to the subscheme $V_t\subseteq X$ defined by $\a^t$,
to $\mu_t=\mu t$, and to $b_t=bt$, we get
$$l(\O_{X,Z}/\a^t)\geq\frac{n^n\mu^{n-1}(\mu+b)}{n!}t^n.$$
Dividing by $t^n$ and passing to the
limit as $t \to \infty$ gives (\ref{e(a)}).

In order to prove~(\ref{l(a)}), we proceed as in~\cite{DEM}.
Passing to the completion, we obtain an ideal $\^\a$ in 
$\^ \O_{X,Z}$. We identify $\^ \O_{X,Z}$ with $K[[x_1,\dots,x_n]]$
via a fixed isomorphism, where $K$ is the residue field of $\O_{X,Z}$.
Moreover, we may choose the local
coordinates so that the image of an equation $h$ defining $H$ in $\O_{X,Z}$
is $x_n$. Since $\^\a$ is zero dimensional,
we can find an ideal $\bb \subset R = K[x_1,\dots,x_n]$,
which defines a scheme supported at the origin, and 
such that $\^\bb = \^\a$. 

If $V'$, $H'\subset{\mathbb A}^n$ are defined by $\bb$ and $x_n$, respectively,
then $({\mathbb A}^n,\frac{1}{\mu}(V'-b\. H'))$ 
is not log terminal at the origin. 
We write
$\mu =  r/s$,
for some $r,s \in \N$, and we may clearly assume that $sb\in\N$.
Consider the ring $S = K[x_1,\dots,x_{n-1},y]$, and the inclusion
$R\subseteq S$ which takes $x_n$ to $y^r$.
This determines a cyclic covering of degree $r$
$$
M := \Spec S \to N := {\mathbb A}^n=\Spec R,
$$
with ramification divisor defined by $(y^{r-1})$.

For any ideal $\cc\subset R$, we put $\tilde\cc:=\cc S$. 
If $W$ is the scheme defined by $\cc$, then we denote by $\widetilde{W}$
the scheme defined by $\tilde\cc$.
In particular,
if $H''\subset M$ is defined by $(y)$, then $\widetilde{H'}=rH''$.
It follows from \cite{ein1}, Proposition~2.8
(see also \cite{Laz}, Section 9.5.E) that 
$(N,\frac{1}{\mu}(V'-b\.H'))$ is not log terminal at the origin in $N$ 
if and only if
$(M,\frac{1}{\mu}\cdot\widetilde{V'}-(sb+r-1)H'')$ is not log terminal
at the origin in $M$.

We write the rest of the proof in 
the language of multiplier ideals, for which we refer
to \cite{Laz}. We use the formal exponential notation for these ideals.
If $\tilde\bb$ is the ideal defining $\widetilde{V'}$, then the above
non log terminality condition on $M$ can be interpreted as saying that
\begin{equation}\label{J}
y^{bs+r-1} \not \in \J(\tilde \bb^{1/\mu}).
\end{equation}

We choose a monomial order in $S$,
with the property that
$$
x_1 > \dots > x_{n-1} > y^{bs+r-1}.
$$
This induces flat deformations to monomial ideals (see \cite{Eisenbud}, 
Chapter 15). 
For an ideal $\dd \subseteq S$, we write the degeneration as
$\dd_t \to \dd_0$, where $\dd_t \cong \dd$ for $t \ne 0$ and 
$\dd_0 =: \ini(\dd)$ is a monomial ideal.

We claim that
\begin{equation}\label{in(J)}
y^{bs+r-1} \not \in \ini(\J(\tilde \bb^{1/\mu})).
\end{equation}
Indeed, suppose that $y^{bs+r-1} \in \ini(\J(\tilde \bb^{1/\mu}))$. Then we
can find an element $f \in \J(\tilde \bb^{1/\mu})$ such that 
$\ini(f) = y^{bs+r-1}$.
Because of the particular monomial order we have chosen, 
$f$ must be a polynomial in $y$ of degree $bs+r-1$. On the other hand,
$\J(\tilde \bb^{1/\mu})$ defines a scheme which is supported at the origin
(or empty), since so does $\tilde \bb$.
We deduce that $y^i\in\J(\tilde\bb^{1/\mu})$, 
for some $i\leq bs+r-1$, which contradicts (\ref{J}).

\begin{lem}\label{in(J(c))vJ(in(c))}
For every ideal $\dd \subseteq S$, and every $c\in{\mathbb Q}_+^*$, we have
$$
\ini(\J(\dd^c)) \supseteq \J(\ini(\dd)^c).
$$
\end{lem}

\begin{proof}
Consider the family $\pi : \MM = \A^n \times T \to T$,
with $T = \A^1$, and the ideal 
$\DDD\subset\O_{\MM}$ corresponding to the degeneration of 
$\dd$ described above.
If $U$ is the complement of the origin in $T$, then
there is an isomorphism 
$$(\pi^{-1}(U), \DDD\vert_{\pi^{-1}(U)})\simeq ({\mathbb A}^n\times U,
{\rm pr}_1^{-1}\dd).$$

Via this isomorphism we have $\J(\pi^{-1}(U),\DDD^c)
\simeq {\rm pr}_1^{-1}(\J(\dd^c))$.
Since the family degenerating to the initial ideal is flat, we deduce easily that
$$\J(\MM,\DDD^c)\cdot \O_{\pi^{-1}(0)}\subseteq\ini(\J(\dd^c)).$$
On the other hand, the Restriction Theorem (see \cite{Laz}) gives
$$\J(\ini(\dd)^c)=\J((\DDD\vert_{\pi^{-1}(0)})^c)\subseteq\J(\MM,\DDD^c)\cdot
\O_{\pi^{-1}(0)}.$$
If we put together the above inclusions, we get the assertion of the lemma.
\end{proof}

Note that the monomial order on $S$ induces
a monomial order on $R$, and that $\widetilde{\ini(\bb)}=\ini(\tilde\bb)$.
Indeed, the inclusion $\widetilde{\ini(\bb)}\subseteq\ini(\tilde\bb)$
is obvious, and the corresponding subschemes have the same length
$r\cdot l(R/\bb)$.

On the other hand, Lemma~\ref{in(J(c))vJ(in(c))} and~(\ref{in(J)}) give
$$
y^{bs+r-1} \not \in \J(\ini(\tilde \bb)^{1/\mu}).
$$
Applying again Proposition 2.8 in \cite{ein1}, in the other direction,
 takes us back in $R$: we deduce that $(N, \frac{1}{\mu}(W-b\cdot H'))$
is not log terminal at the origin, where $W\subset N$ is defined by $\ini(\bb)$.
Since $l(\O_{X,Z}/\a)=l(R/\bb)=l(R/\ini(\bb))$, we have reduced the proof
of (\ref{l(a)}) to the case when $\a$ is a monomial ideal.
In this case, we have in fact a stronger statement, which we prove
in the lemma below; therefore the proof of
Theorem~\ref{l(a)-e(a)} is complete.
\end{proof}

The following is the natural generalization of Lemma~2.1 in \cite{DEM}.

\begin{lem}\label{monomial}
Let $\a$ be a zero dimensional
monomial ideal in the ring $R = K[x_1,\dots,x_n]$, defining a scheme $V$.
Let $H_i$ be the hyperplane defined by $x_i=0$. 
We consider $\mu\in{\mathbb Q}_+^*$ and $b_i\in{\mathbb Q}$, such that
$\mu\geq\max_i\{b_i\}$.
If the pair $({\mathbb A}^n,
\frac{1}{\mu}(V+\sum_ib_iH_i))$ is not log terminal, then
$$
l(R/\a)\geq\frac{n^n}{n!} \. \prod_{i=1}^n (\m - b_i),
$$
and the inequality is strict if $n\geq 2$.
\end{lem}

\begin{proof}
We use the result in \cite{ELM} which gives the
condition for a monomial pair, with possibly negative coefficients,
to be log terminal. This generalizes the formula for the log canonical
threshold of a monomial ideal from \cite{Ho}.
It follows from \cite{ELM} that 
$(X,\frac{1}{\mu}(V+\sum_ib_iH_i))$ is not log terminal if and only if
there is a facet of the Newton polytope associated to $\a$ such that, 
if $\sum_i u_i/a_i = 1$ is the equation of the hyperplane 
supporting it, then 
$$
\sum_{i=1}^n \frac{\m- b_i}{a_i}\leq 1.
$$
Applying the inequality between the arithmetic mean and the geometric mean
of the set of nonnegative numbers $\{(\m - b_i)/a_i\}_i$, we deduce
$$
\prod_i a_i \geq n^n \. \prod_i (\m - b_i).
$$
We conclude using the fact that 
$n! \. l(R/\a)\geq \prod_i a_i$,
and the inequality is strict if $n\geq 2$
 (see, for instance, Lemma~1.3 in~\cite{DEM}).
\end{proof}

\bigskip

\section{Log canonical thresholds of affine cones}

In this section we give a lower bound for the log canonical threshold
of a subscheme $V\subset\A^n$, cut out by homogeneous 
equations of the same degree.
The bound involves the dimension of the non log terminal locus
of $(\A^n,c\cdot V)$, 
where $c=\lc(\A^n,V)$. Moreover, we characterize the 
case when we have equality.
In the particular case when $V$ is the affine cone over a projective 
hypersurface
with isolated singularities, this proves a conjecture of Cheltsov and Park
from~\cite{CP}.

The main ingredient we use for this bound is a formula for the log
canonical threshold in terms of jet schemes, from \cite{Mu2}.
Recall that for an arbitrary scheme $W$, 
of finite type over the ground field $k$, the $m$th jet scheme $W_m$
is again a scheme of finite type over $k$ characterized by
$${\rm Hom}({\rm Spec}\,A, W_m)\simeq{\rm Hom}({\rm Spec}\,A[t]/(t^{m+1}),
W),$$
for every $k$-algebra $A$. Note that $W_m(k)=
{\rm Hom}({\rm Spec}\,k[t]/(t^{m+1}), W),$ and in fact, we will be
interested only in the dimensions of these spaces.
For the basic properties of the jet schemes, we refer to
\cite{Mu1} and \cite{Mu2}.

\begin{thm}\label{ingred}{\rm (\cite{Mu2}, 3.4)}
If $X$ is a smooth, connected variety of dimension $n$, and 
if $V\subset X$
is a subscheme, then the log canonical threshold of $(X,V)$
is given by
$$\lc(X,V)=n-\sup_{m\in\N}\frac{\dim\,V_m}{m+1}.$$
Moreover, there is $p\in\N$, depending on the numerical data given
by a log resolution of $(X,V)$, such that $\lc(X,V)=n-(\dim\,V_m)/(m+1)$
whenever $p\mid (m+1)$.
\end{thm}

\smallskip

For every $W$ and every $m\geq 1$, there are canonical projections
$\phi^W_m:W_m\longrightarrow W_{m-1}$ induced by the truncation
homomorphisms $k[t]/(t^{m+1})\longrightarrow k[t]/(t^m)$.
By composing these projections we get morphisms $\pi^W_m:W_m
\longrightarrow W$. When there is no danger of confusion, 
we simply write $\phi_m$ and $\pi_m$.

If $W$ is a smooth, connected variety, then $W_m$ is smooth, connected,
and $\dim\,W_m=(m+1)\dim\,W$, for all $m$. It follows from definition
that taking jet schemes commutes with open immersions. In particular,
if $W$ has pure dimension $n$, then $\pi_m^{-1}(W_{\rm reg})$
is smooth, of pure dimension $(m+1)n$.

Recall that the non log terminal locus of a pair is the union of all centers
of non log terminality. 
In other words, its complement is the largest open subset
over which the pair is log terminal.
Theorem~\ref{ingred} easily gives a description via jet schemes
of the non log terminal locus of a pair which is log canonical, but is not
log terminal.
Suppose that $(X, V)$ is as in the theorem, and
let $c=\lc(X,V)$. We say that an irreducible component $T$ of $V_m$
(for some $m$) computes $\lc(X,V)$ if $\dim(T)=(m+1)(n-c)$. Note that basic
results on jet schemes show that for every irreducible component $T$ of
$V_m$, the projection $\pi_m(T)$ is closed in $V$ (see \cite{Mu1}).
It follows from Theorem~\ref{ingred}
that if $W$ is an irreducible
component of $V_m$ that computes the log canonical threshold of $(X, V)$
then $\pi_m(W)$ is contained in the 
non log terminal locus of $(X,c\cdot V)$
(see also \cite{ELM}).

\medskip

For future reference, we record here two lemmas.
For $x \in \R$, we denote by $[x]$ the largest integer $p$
such that $p \le x$.

\begin{lem}\label{fiber}{\rm (\cite{Mu1}, 3.7)}
If $X$ is a smooth, connected variety of dimension $n$, $D\subset X$
is an effective divisor, and $x\in D$ is a point with $e_xD=q$,
then 
$$\dim(\pi^D_m)^{-1}(x)\leq mn-[m/q],$$
for every $m\in\N$.
\end{lem}

In fact, the only assertion we will need from Lemma~\ref{fiber}
is that $\dim\,(\pi^D_m)^{-1}(x)\leq mn-1$, if $m\geq q$, which follows easily
from the equations describing the jet schemes (see \cite{Mu1}).

\smallskip

\begin{lem}\label{semicont}{\rm (\cite{Mu2} 2.3)}
Let $\Phi : {\mathcal W}\longrightarrow S$ be a family of schemes,
and let us denote the fiber $\Phi^{-1}(s)$ by ${\mathcal W}_s$.
If $\tau:S\longrightarrow{\mathcal W}$ is a section of $\Phi$, then 
the function
$$f(s)=\dim(\pi_m^{{\mathcal W}_s})^{-1}(\tau(s))$$
is upper semi-continuous on the set of closed points of $S$, for every $m\in\N$.
\end{lem}

\bigskip

The following are the main results in this section.

\begin{thm}\label{lower_bound1}
Let $V\subset\A^n$ be a subscheme whose ideal is generated by homogeneous
polynomials of degree $d$. Let $c=\lc(\A^n,V)$, and let $Z$ be the 
non log terminal locus of $(\A^n, c\cdot V)$. If $e=\codim(Z,\A^n)$, then 
$c\geq e/d$.
\end{thm}

\begin{thm}\label{equality_case1}
With the notation in the previous theorem, $c=e/d$ if and only if
$V$ satisfies the following three properties:
\begin{enumerate}
\item
$Z = L$ is a linear subspace of codimension $e$. 
\item
$V$ is the pull back of a closed subscheme $V'\subset\A^n/L$, 
which is defined
by homogeneous polynomials of degree $d$ and such that
$\lc(\A^n/L, V') =e/d$.
\item
The non log terminal locus
of $(\A^n/L, e/d\cdot V')$ is just the origin. 
\end{enumerate}
\end{thm}

\begin{proof}[Proof of Theorem~\ref{lower_bound1}]
If
$\pi_m:V_m\longrightarrow V$ is the canonical projection,
then we have an isomorphism
\begin{equation}\label{isom}
\pi_m^{-1}(0)\simeq V_{m-d}\times\A^{n(d-1)},
\end{equation}
for every $m\geq d-1$ (we put $V_{-1}=\{0\}$).
Indeed, for a $k$-algebra $A$, an $A$-valued point
of $\pi_m^{-1}(0)$ is a ring homomorphism
$$\phi:k[X_1,\ldots,X_n]/(F_1,\ldots,F_s)\longrightarrow
A[t]/(t^{m+1}),$$
such that $\phi(X_i)\in(t)$ for all $i$.
Here $F_1,\ldots, F_s$ are homogeneous equations 
of degree $d$, defining $V$. Therefore we can write
$\phi(X_i)=tf_i$, and $\phi$ is a homomorphism if and only if
the classes of $f_i$ in $A[t]/(t^{m+1-d})$ define an
$A$-valued point of $V_{m-d}$. But $\phi$ is uniquely 
determined by the classes of $f_i$ in $A[t]/(t^m)$, so this proves 
the isomorphism in equation~(\ref{isom}).

By Theorem~\ref{ingred}, we can find $p$ such that
$$\dim\,V_{pd-1}=pd(n-c).$$
Let $W$ be an irreducible component of $V_{pd-1}$
computing $\lc(X,V)$, so $\dim\,W=pd(n-c)$ and $\pi_{pd-1}(W) \subset
Z$. 
By our hypothesis,
$\dim\pi_{pd-1}(W)\leq n-e$. Therefore Lemma~\ref{semicont} gives
\begin{equation}\label{inequality1}
pd(n-c)=\dim\,W\leq\dim\pi_{pd-1}^{-1}(0)+n-e=
\dim V_{(p-1)d-1}+(d-1)n+n-e,
\end{equation}
where the last equality follows from (\ref{isom}).
Another application of Theorem~\ref{ingred} gives
\begin{equation}\label{inequality2}
\dim\,V_{(p-1)d-1}\leq (p-1)d(n-c).
\end{equation}
Using this and (\ref{inequality1}), we get
$c\geq e/d$.
\end{proof}

\begin{proof}[Proof of Theorem~\ref{equality_case1}]
We use the notation in the above proof.
Since $c=e/d$, we see that
in both equations~(\ref{inequality1}) and (\ref{inequality2})
we have, in fact, equalities. The equality in (\ref{inequality2}) shows
that $\dim V_{(p-1)d-1}=(p-1)d(n-c)$, so
we may run the same argument with 
$p$ replaced by $p-1$.
Continuing in this way, we see that we may suppose that $p=1$.
In this case, the equality in (\ref{inequality1}) shows that
for some irreducible component $W$ of $V_{d-1}$, with $\dim W=dn-e$,
we have $\dim \pi_{d-1}(W)=n-e$. It follows that if $Z_1:=\pi_{d-1}(W)$,
then $Z_1$ is an irreducible component of $Z$.

Fix $x\in Z_1$. If ${\rm mult}_xF\leq d-1$,
for some degree $d$ polynomial $F$ in the ideal of $V$, then
Lemma~\ref{fiber} would give $\dim\,\pi_{d-1}^{-1}(x)\leq (d-1)n-1$.
This would imply $\dim\,W\leq n-e+(d-1)n-1$,
a contradiction. Therefore we must have ${\rm mult}_xF\geq d$,
for every such $F$. 

Recall that we have degree $d$ generators of the ideal of
$V$, denoted by $F_1,\ldots,F_s$. 
Let $L_i = \{x \in \A^n | {\rm mult}_xF_i =d\}$, for $i\leq s$. 
By the B\'{e}zout
theorem, $L_i$ is a linear space. If $L= \bigcap_{i=1}^s L_i$, then 
$Z_1 \subset L$. On the other hand, 
by blowing-up along $L$, we see that $L$ is contained
in the non log terminal locus of $(\A^n, c\cdot V)$. Therefore
$Z_1 =L$. Let $z_1, ..., z_e$ be the linear forms defining $L$. Then
each $F_i$ is a homogeneous polynomial of degree $d$ 
in $z_1, ..., z_e$. This shows
that $V$ is the pull back of a closed subscheme $V'\subset\A^n/L$,
defined by $F_1,..., F_s$. Since the projection map $\pi: \A^n
\longrightarrow \A^n/L$ is smooth and surjective, we see that 
$\lc(\A^n/L, V') = \lc(\A^n, V)$ and that the non log terminal locus 
of $(\A^n, \frac{e}{d}\cdot V)$ is just the pull-back of the corresponding locus
for the pair $(\A^n/L, e/d\cdot V')$.
Note that the non log terminal locus of $(\A^n/L, e/d\cdot V)$
is defined by an homogeneous ideal. By dimension considerations, 
we conclude that this locus consists just of the origin, so
$Z= L$. 

Conversely, if $V$ is the pull back of a closed subscheme from
$\A^n/L$ as described in the theorem, one checks that 
$\lc(\A^n, V) = e/d$ and that the corresponding
non log terminal locus is just $L$.
\end{proof}

Let $V'$ be a closed subscheme of $\P^{n-1}$ defined by degree $d$ 
homogeneous polynomials $F_1,\ldots, F_s$, and let $V$ be the closed
subscheme in $\A^n$ defined by the same set of
polynomials. Let $c =\lc(\P^{n-1}, V')$, and let $Z'$ be the 
non log terminal locus of
$(\P^{n-1}, c\cdot V')$. Suppose that the codimension of $Z'$
in $\P^{n-1}$ is $e$.

\begin{cor}\label{proj_case}
With the above notation, $\lc(\P^{n-1}, V') \ge e/d$. Moreover,
if we have equality, then $V'$ is the cone over a scheme in some $\P^{e-1}$.
\end{cor}

\begin{proof}
Note that 
$$\lc(\P^{n-1}, V') = \lc(\A^n-\{0\}, V-\{0\}) \ge
\lc(\A^n, V).$$
Now the first assertion follows from Theorem~\ref{lower_bound1}.

If $\lc(\P^{n-1}, V') = e/d$, then $\lc(\A^n, V) = e/d$ and
the non log terminal locus of $(\A^n, \frac{e}{d}\cdot V)$ is a
linear space $L$ of codimension 
$e$. If $z_1,..., z_e$ are the linear forms defining $L$, then each
$F_i$ is a homogeneous polynomial of degree $d$ in $z_1, ..., z_e$.
Therefore $V'$ is the cone with center $L$
over the closed subscheme of $\P^{e-1}$
defined by $F_1,\ldots,F_s$. 
\end{proof}

\begin{rmk}
In \cite{CP}, Cheltsov and Park studied the log canonical threshold
of singular hyperplane sections of smooth, projective hypersurfaces. 
If $X\subset{\mathbb P}^n$
 is a smooth hypersurface of degree $d$, and if ${V}
=X\cap H$, for a hyperplane $H$, then 
they have shown that
\begin{equation}\label{ineq_CP}
\lc(X,{V})\geq\min\{(n-1)/d,1\}.
\end{equation} 
It follows from Theorem~\ref{ingred} that $\lc(X, V) = \lc 
(\P^{n-1}, V)$. As
it is well known that 
${V}$ has isolated singularities, if we apply 
the first assertion in Corollary~\ref{proj_case},
then we recover the result in \cite{CP}. 

Cheltsov and Park have conjectured in their setting that if $d\geq n$, then
equality holds in
(\ref{ineq_CP}) if and only if ${V}$ is a cone. They have shown that
their conjecture would follow from the Log Minimal Model Program.
The second assertion in Corollary~\ref{proj_case}
proves, in particular, their conjecture.
\end{rmk}

\bigskip

\section{Application to birational rigidity}

Using the bounds on log canonical thresholds
from the previous
sections, we prove now the birational rigidity of 
certain Fano hypersurfaces.
We recall that a Mori fiber space $X$ is
called {\it birationally superrigid} if any 
birational map $\f : X \rat X'$ to another Mori fiber space $X'$ 
is an isomorphism. For the 
definition of Mori fiber space and for another notion
of rigidity, we refer to~\cite{Co2}. Note that Fano manifolds 
having N\'{e}ron-Severi group of rank 1 are
trivially Mori fiber spaces.
Birational superrigidity is a very strong condition:
it implies that $X$ is not rational,
and that $\Bir(X) = \Aut(X)$.  Note that if $X$ is a smooth hypersurface
of degree $N$ in $\P^N$ ($N \ge 4$), then $X$ has no nonzero 
vector fields. Therefore if $X$ is birationally superrigid, 
then the birational invariant $\Bir(X)$ is a finite
group.

The following theorem is the main result of this section.

\begin{thm}\label{X_N}
For any integer $4 \le N \le 12$,
every smooth hypersurface $X = X_N \subset \P^N$ of degree $N$
is birationally superrigid. 
\end{thm}

The case $N=4$ of the above theorem is due to Iskovskikh and Manin
(see~\cite{IM}). The case $N=5$ was proven by Pukhlikov
in~\cite{Pu2}, while the cases $N=6,7,8$ were established by 
Cheltsov in~\cite{Ch2}.
Birational superrigidity of smooth hypersurfaces of degree $N$
in $\P^N$ (for $N \ge 5$) was conjectured by Pukhlikov in~\cite{Pu5}, where
the result is established under a suitable condition
of regularity on the equation defining the hypersurface. 
We remark that there is an attempt 
due to Pukhlikov in \cite{Pu1} to prove the general case (for $N \ge 6$).
Despite a gap in the proof
(see the remark below), we believe that the method therein
could lead in the future to the result.
In fact, the proof given below for Theorem~\ref{X_N}
follows his method, and our contribution is mainly in simplifying
and solidifying his argument.

\begin{rmk}\label{gap}
The following gives a counterexample to Corollary~2 in \cite{Pu1}.
Let $Q\subset{\mathbb P}^4$ be a cone over a twisted cubic,
and let $\pi_a: Q\longrightarrow R=\pi_a(Q)$ be the projection from
an arbitrary point
$a\in{\mathbb P}^4\setminus Q$; 
note that $R$ is the cone over a singular plane cubic.
 If $p$ is the vertex of $Q$, then the restriction of
$\pi_a$ to any punctured neighbourhood of $p$ in $Q$ can not preserve
multiplicities, as $q=\pi_a(p)$ lies on a one dimensional component
of the singular locus of $R$.
\end{rmk}

Before proving the above theorem, we recall the following 
result, due to Pukhlikov:

\begin{prop}\label{pu1}{\rm (\cite{Pu1}, Proposition~5)}
Let $X \subset \P^N$ be a smooth hypersurface, and let
$Z$ be an effective 
cycle on $X$, of pure codimension $k < \frac 12 \dim X$.
If $m \in \N$ is such that 
$Z \equiv m \.c_1(\O_X(1))^k$,  
then $\dim \{ x \in Z \mid e_xZ > m \} < k$.
\end{prop}

\begin{rmk}\label{pu1_rmk}
Because we have assumed
$k<\frac 12 \dim X$, the existence of $m$ as in the proposition
follows from Lefschetz Theorem.
One can check that the proof of Proposition~\ref{pu1}
extends to the case $k = \frac 12 \dim X$, if we assume
that such $m$ exists.
Note also that the statement is trivially true
if $k > \frac 12 \dim X$.
\end{rmk} 

We need first a few basic properties which allow us
to control multiplicities when restricting to general hyperplane
sections, and when projecting to lower dimensional linear subspaces. The
following proposition must be well known, but we include a proof for the 
convenience of the readers. We learned this proof, 
which simplifies our original arguments, from Steve Kleiman.

\begin{prop}\label{int_mult}
Let $Z\subset\P^n$ be an irreducible projective variety. 
If $H \subset Z$ is a general hyperplane section, then
$e_pH = e_pZ$ for every $p \in H$.
\end{prop}

\begin{proof}
As observed by Whitney (e.g., see~\cite{Kl}, page 219), 
at any point $p \in Z$, the fiber over $p$ of the conormal variety of $Z$,
viewed as a linear subspace of $(\P^n)^*$,
contains the dual variety of every component of the
embedded projective tangent
cone $C_pZ$ of $Z$ at $p$. A hyperplane section $H$ of $Z$
satisfies $e_pH = e_pZ$ if the hyperplane meets $C_pZ$
properly. Therefore, this equality holds for every point
$p$ in $H$ whenever $H$ is cut out by a hyperplane
not in the dual variety of $Z$. 
\end{proof}

In the next two propositions, we consider a (possibly reducible)
subvariety  $Z \subset \P^{n+s}$, of pure dimension $n-1$,
for some $n \ge 2$ and $s\geq 1$, and take
a general linear projection $\p : \P^{n+s} \setminus \LL \to \P^n$. 
Here $\LL$ denotes the center of projection, that is an 
$(s-1)$ dimensional linear space.
We put $T = \p(Z)$ and $g = \p|_Z : Z \to T$. 
It is easy to see that since $\LL$ is general, $g$ is a finite
birational map.
For convenience, we put $\dim(\emptyset)=-1$.

\begin{prop}\label{proj_mult1}
With the above notation, consider the set
$$
\D = \Big\{q \in T \mid e_q T > 
\sum_{p \in g^{-1}(q)} e_p Z \Big\}. 
$$
If the projection is chosen with suitable generality,
then $\codim(\D,{\mathbb P}^n) \ge 3$.
\end{prop} 

\begin{proof}
Note that $e_qT \ge \sum e_pZ$ for every $q \in T$, 
the sum being taken over all points $p$ over $q$. Moreover,
for a generic projection, every irreducible component of $Z$
is mapped to a distinct component of $T$. 
Therefore, by the linearity of the multiplicity, we may assume that
$Z$ is irreducible.

Let $\D' \subset T$ be the set of points $q$,
such that for some $p$ over $q$, the intersection of
the $s$ dimensional linear space $\ov{\LL q}$ with the
embedded projective tangent cone $C_pZ$
of $Z$ at $p$, is at least one dimensional. We claim that 
$\codim (\D',{\mathbb P}^n) \geq 3$.
Indeed, it follows from the theorem on generic flatness that
there is a  stratification $Z=Z_1 \sqcup \dots \sqcup Z_t$ 
by locally closed subsets
such that, for every $1 \le j \le t$, the incidence set
$$
I_j = \{ (p,x) \in Z_j\times\P^{n+s} \mid x \in C_pZ \}
$$
is a (possibly reducible) quasi-projective variety of dimension 
no more than $2 \dim Z= 2n-2$. 
Let $\pr_1$ and $\pr_2$ denote the projections of $I_j$ to the first
and to the second factor, respectively.
It is clear that the set of those $y\in{\mathbb P}^{n+s}$, with
$\dim \pr_2^{-1}(y)=\tau$ has dimension at most $\max\{2n-2-\tau,-1\}$,
for every $\tau\in{\mathbb N}$. Since
$\LL$ is a general linear subspace of dimension $s-1$, 
it intersects a given $d$ dimensional closed subset in a set of
dimension $\max\{d-n-1,-1\}$. 
Hence $\dim \pr_2^{-1}(\LL)\leq n-3$, and
therefore 
$\dim(\pr_1(\pr_2^{-1}(\LL))) \le n-3$.
As this is true for every $j$, we deduce $\codim(\D',{\mathbb
  P}^n)\geq 3$.
Thus, in order to prove the proposition, it is enough to show that
$\D \subseteq \D'$. 
 
For a given point $p \in Z$, let $L_p \subset \P^{n+s}$ be an 
$(s+1)$ dimensional linear subspace passing through $p$.
Let $\mm_p$ be the maximal ideal of $\O_{Z,p}$, and let
$\PP \subset \O_{Z,p}$ be the ideal locally defining $L_p \cap Z$. 
If $L_p$ meets the tangent cone $C_pZ$ of $Z$ at $p$ properly,
then the linear forms defining $L_p$
generate the ideal of the exceptional divisor of the 
blow up of $Z$ at $p$. Therefore $e(\mm_p)=e(\PP)$.

Consider now some $q \in T \setminus \D'$. Let $L_q \subset \P^n$ 
be a general line passing through $q$,
and let $\QQ \subset \O_{T,q}$ be the ideal generated by
the linear forms vanishing along $L_q$.
We denote by $L$ the closure of $\p^{-1}(L_q)$ in $\P^{n+q}$.
For every $p \in g^{-1}(q)$, let $\PP \subset \O_{Z,p}$ be the ideal 
generated by the linear forms vanishing along $L$. 
Since $L_q$ is general and $q \not \in \D'$,
we may assume that $L$ intersects $C_pZ$ properly,
hence $e(\mm_p) = e(\PP)$.
On the other hand, if $\mm_q$ is the maximal ideal of $\O_{T,q}$,
then $\QQ\subseteq\mm_q$, which gives 
$$
\PP=\QQ\cdot\O_{Z,p} \subseteq \mm_q\cdot\O_{Z,p}\subseteq\mm_p. 
$$
Therefore 
$e(\mm_p)=e(\mm_q\cdot\O_{Z,p})$ for every $p$ as above,
hence $q\not\in\D$,
by \cite{Fulton}, Example~4.3.6.
\end{proof}

\begin{prop}\label{proj_mult2}
With the notation  in Proposition~\ref{proj_mult1}, 
consider the set
$$
\S = \S(Z,\p):= \{q \in T \mid \text{$g^{-1}(q)$ has al least 3
distinct points} \}.
$$
If the projection is sufficiently general, then  
$\codim(\S,{\mathbb P}^n) \geq 3$. 
\end{prop} 

\begin{proof}
We have $\codim(\S,{\mathbb P}^n)\geq 3$ 
if and only if $\S \cap P = \emptyset$
for every general plane
$P\subset\P^n$. Pick one general plane $P$,
let $P' \;(\cong \P^{s+2})$ be the closure of $\p^{-1}(P)$ in $\P^{n+s}$,
and let $\p'$ be the restriction of $\p$ to $P' \setminus \LL$.
If $Z' = Z \cap P'$, then 
$Z'$ is a (possibly reducible) curve, and its multisecant variety
is at most two dimensional (see, for example, \cite{FOV}, Corollary~4.6.17).
Note that $\LL$ is general in $P'$. Indeed, choosing 
the center of projection $\LL$ general
in $\P^{n+s}$, and then picking $P$ general in $\P^n$ is equivalent to first 
fixing a general $(s+2)$-plane $P'$ in $\P^{n+s}$ and then choosing $\LL$ 
general in $P'$.
Therefore we conclude that $\S \cap P$, which is the same as 
$\S(Z',\p')$, is empty.
\end{proof}

\begin{proof}[Proof of Theorem~\ref{X_N}]
By adjunction, $\O_X(-K_X)\simeq\O_X(1)$.
Let $\f : X \rat X'$ be a birational map from $X$ to a Mori fiber space $X'$,
and assume that $\f$ is not an isomorphism.
By the Noether-Fano inequality
(see~\cite{Co1} and~\cite{Is}, or~\cite{Ma}), we find a linear
subsystem $\H \subset |\O_X(r)|$, with $r \geq 1$, whose base scheme
$B$ has codimension $\ge 2$, and such that the pair $(X,\frac 1r\.B)$ is not
canonical. We choose $c<\frac{1}{r}$, such that $(X,c\cdot B)$ is still
not canonical, and
let $C \subset X$ be a center of 
non canonicity for $(X,c \.B)$.
Note that $C$ is a center of non canonicity also for 
the pairs
$(X,c \. D)$ and $(X,c \. V)$, 
where $V = D \cap D'$ and $D$, $D' \in \H$ are
two general members. Applying Proposition~\ref{pu1} for $Z=D$ and $k=1$,
we see that the multiplicity of $D$ is $\leq r$ on an open subset
whose complement has dimension zero. 
On this open subset $(X,c\.D)$
is canonical 
(see, for example, \cite{Ko} 3.14.1). 
Therefore $C = p$, a point of $X$.

Let $Y$ be a general hyperplane section of $X$ containing $p$.
Then $p$ is a center of non log canonicity for $(Y,c \. B|_Y)$.
Note that $Y$ is a smooth hypersurface of degree $N$
in $\P^{N-1}$.  
Let $\p : \P^{N-1} \setminus \LL \to \P^{N-3}$ be a general
linear projection, where the center of projection $\LL$ is a line. 
We can assume that
the restriction of $\p$ to 
each irreducible component of $V|_Y$ is finite and birational. 
Note that $\p_*[V|_Y]$ is a divisor in $\P^{N-3}$ of degree $Nr^2$.
If $\tilde  Y = \Bl_{\LL \cap Y} Y$, then we get a morphism
$f : \tilde  Y \to \P^{N-3}$.
If we choose $\LL$ general enough, 
then we can find an open set $U \subset \P^{N-3}$,
containing the image $q$ of $p$, 
such that $f$ restricts to a smooth (proper) morphism
$f^{-1}(U) \to U$. Applying Theorem~\ref{thm1}, we deduce 
that the pair
\begin{equation}\label{pair}
\(\P^{N-3}, \frac{c^2}4 \. \p_*[V|_Y] \)
\end{equation}
is not log terminal at $q$.

We claim that
\begin{equation}\label{dim_bound}
\dim \{y \in \pi(V|_Y) \mid e_y(\p_*[V|_Y]) > 2r^2 \}
\le \max \{ N-6,0 \}.
\end{equation}
Indeed, by Propositions~\ref{proj_mult2} and~\ref{proj_mult1}, the map 
$\Supp([V|_Y]) \to \Supp(\p_*[V|_Y])$ is at most 2 to 1 and preserves
multiplicities outside a set, say $\D \cup \S$, 
of dimension $\le \max\{N-6,-1\}$.
This implies that,
for each $y$ outside the set $\D \cup \S$, 
$e_y(\p_*[V|_Y]) = \sum e_x([V|_Y])$, 
where the sum is taken over the points $x$ over $y$, and this 
sum involves at most two non-zero terms.
Then~(\ref{dim_bound}) follows from the fact that,
by Propositions~\ref{pu1} and ~\ref{int_mult}
(see also Remark~\ref{pu1_rmk}),
the set of points $x$ for which
$e_x[V|_Y] > r^2$ is at most zero dimensional.

Note that the pair~(\ref{pair}) is log terminal at every point
$y$ where $e_y(\p_*[V|_Y]) \le 4r^2$. If $4 \le N \le 6$, 
we deduce that the pair is log terminal outside a 
zero dimensional closed subset.
In this case, Corollary~\ref{proj_case} gives
$c^2/4 \ge (N-3)/(Nr^2)$. Since $c < 1/r$, this implies $N < 4$, 
a contradiction.
If $7 \le N \le 12$, then we can only conclude that 
the pair~(\ref{pair}) is log terminal outside a closed subset
of codimension at least $3$.
This time the same corollary gives
$c^2/4 \ge 3/(Nr^2)$, which implies $N > 12$. 
This again contradicts our assumptions,
so the proof is complete.
\end{proof}

\providecommand{\bysame}{\leavevmode \hbox \o3em
{\hrulefill}\thinspace}

\end{document}